\documentclass[12pt]{article} 
\usepackage{amsmath}
\usepackage[dvips]{graphicx}
\usepackage{amsfonts}
\usepackage{amsopn}
\usepackage{amsthm}
\def\bb{\bf}

\def\refl{\cal{R}}

\def\xy{\langle x, y \rangle}
\def\XY{\langle X, Y \rangle}

\def\Reals{{\bb R}}
\def\Realsb{{\bar{\Reals}}}

\def\Rn{{\Reals}^{n+1}}

\def\Sn{ {\bf S}^n} 

\def\calO{{\cal O}}

\newcommand{\p}{\rho}

\newcommand{\g}{\gamma}

\renewcommand{\a}{\alpha}

\newcommand{\ra}{\rangle}

\newcommand{\la}{\langle}

\newtheorem{theorem}{Theorem}
\newtheorem{definition}[theorem]{Definition}

\newtheorem{lemma}[theorem]{Lemma}
\newtheorem{proposition}[theorem]{Proposition}
\newtheorem{remark}[theorem]{Remark}
\title {A Minkowski-style theorem for focal functions of
compact convex reflectors}
\author{Vladimir I. Oliker\thanks{The research of the author
was partially supported by the National Science Foundation grant DMS-04-05622,
the Air Force Office of Scientific Research under contract FA9550-05-C-0058
and by a grant from the Emory University Research Committee.}\\
Department of Mathematics and Computer Science,\\
 Emory University, Atlanta, Georgia\\
oliker@mathcs.emory.edu}
\date{}
\begin{document}
\maketitle
\pagenumbering{arabic}
\setcounter{section}{0}
\setcounter{subsection}{0}
\begin{abstract} 
This paper continues the study of a class of
 compact convex hypersurfaces
in $\Rn, ~n \geq 1,$ which are boundaries of compact convex 
bodies obtained by taking the intersection of (solid) confocal paraboloids of revolution. Such hypersurfaces are called
reflectors.  In $\Reals^3$ reflectors
arise naturally in geometrical optics and are used in design 
of light reflectors and reflector antennas. They are also important in
rendering problems in computer graphics.
 
The notion of a focal function for reflectors plays a central role similar to that of the Minkowski support function
for convex bodies. In this paper the basic question of when a given function is
a focal function of a convex reflector is answered by establishing necessary and sufficient conditions. In addition, some smoothness properties of reflectors and
of the associated directrix hypersurfaces are also etablished.
\end{abstract}
\section{\protect\bf Introduction}
\label{intro}
A convex reflector is a convex hypersurface which is the boundary of a compact convex body in Euclidean space $\Rn, n\geq 1,$
obtained by taking the intersection of a given set of confocal (solid) paraboloids of revolution. Convex reflectors
arise naturally as solutions to nonlinear second order elliptic partial differential equations of
Monge-Amp\`{e}re type describing conservation laws in geometrical optics \cite{W}, \cite{CO}, \cite{CKO}.
Numerous optics and electromagnetic applications require solutions of such equations; for example, the most
common reflector antennas are designed with the use of optics conservation laws and requre determination
of reflecting surfaces \cite{W}. In addition, it follows from recent results in 
\cite{Caf_alloc:96}, \cite{glimm-oliker-refl2:03}, \cite{glimm-oliker-refl1:03}, \cite{galo:04}, \cite{Wang:inv03} that in variational formulations of such
problems as problems of 
Monge-Kantorovich optimal mass transfer the radial and focal functions of a convex reflector
turned out to be the logarithms of the Kantorovich potentials \cite{KA}. Convex reflectors are also important in inverse problems of rendering in computer graphics when the scattering characteristics of reflectors are prescribed in advance \cite{patow_pueyo}.

These reasons already provide sufficient motivation for a systematic study of this class of convex hypersurfaces.
Moreover, the results in \cite{Hasanis_Koutrou:85}, \cite{refl_geom1} and \cite{CGH:04} show that convex reflectors
are of independent geometric interest and this paper continues the investigations in this direction. 

We describe now briefly the content and organization of the paper. In order to make
the presentation reasonably self-contained, we first review in section \ref{prel} the basic definitions and a few facts about
convex reflectors. More details are provided in \cite{refl_geom1}, but this
paper can be read independently of \cite{refl_geom1}. A new fact, also shown in this section, is that
the support function of a nondegenerate convex reflector is of class $C^1$. 

In section \ref{mink} we study
analytic properties of the $focal$ function of a convex
reflector. In convexity theory the support function is
one of the main tools for studying convex hypersurfaces as any convex 
body can be described by its support function. The classical
theorem of Minkowski asserts that any sublinear function in $\Rn$ is
the support function of a unique convex body \cite{Schneider}. By contrast, a convex reflector
is more conveniently described by the focal function which for each $y$ 
on a unit sphere $\Sn$ defines the focal parameter of a paraboloid of revolution
with axis $y$ supporting to that reflector. 
The main result in section \ref{mink} is Theorem \ref{minkg} which gives
necessary and sufficient conditions for a function on $\Sn$
to be a focal function of a convex reflector. Thus, this theorem can
be viewed as a generalization of the result  of
Minkowski. Our results
also imply that sublinearity by itself (of the appropriately
extended from $\Sn$ to $\Rn$ focal function) is not sufficient
to define a focal function of a convex reflector.

In section \ref{dir} we study the directrix hypersurface of
the reflector which is (up to rescaling) is its pedal hypersurface. In particular,
it is shown that the directrix of a convex reflector is a hypersurface of class $C^1$. A geometric proof of this fact was given earlier in \cite{refl_geom1}; the
proof presented here is analytic and more simple.

\section{\protect\bf Preliminaries}\label{prel}

In $\Rn$ fix a Cartesian coordinate system with 
the origin $\cal{O}$
and
let $\Sn$ be the unit sphere centered at $\cal O$. 
Throughout the paper the term $paraboloid$
means a paraboloid of revolution with focus at
$\cal{O}$ and axis $y \in \Sn$ directed towards its opening. 
Such paraboloid is denoted by 
$P(y)$. No other kind of paraboloids will be considered.
The closed convex subset of $\Rn$ bounded by $P (y)$ and
containing its axis  is denoted by $B (y)$.

The polar radius
of $P(y)$ is given by
\begin{equation}\label{polrad}
\p_y (x) = \frac{\tilde {p}}{1 - \xy}, ~~~x \in \Sn \setminus \{y\},
\end{equation}
where $\tilde {p}$ is a nonnegative real number called 
the {\it focal parameter} of $P (y)$. When $\tilde {p}=0$ it is assumed
that 
the paraboloid reduces to an infinite
ray in direction $y$ emanating from $\calO$. Such paraboloid is called
$degenerate$. When $ \tilde {p}=\infty$ the paraboloid is called $improper$.
In case of an improper paraboloid, $B (y) \equiv \Rn$.

\begin{definition} \label{reflector} Let $\{P(y), ~y \in \Sn\}$ be a family
of confocal paraboloids,
\[ B = \bigcap_{y \in \Sn } B (y),~\mbox{and}~ R = \partial B. \]
If $B$ is compact and contains
interior points (in the topology of $\Rn$), 
the closed convex hypersurface $R$ is called a {\it convex reflector
(with the light source $\cal{O})$}. If $B$ has an empty interior, $R$ is
called degenerate. The set of all nondegenerate convex 
reflectors with the
light source at  $\cal{O}$ is denoted by $\refl$.
\end{definition}

Since a unit vector $y$ and
focal parameter $\tilde {p}$ define a paraboloid $P(y)$ uniquely, it is
convenient to describe  families of confocal paraboloids by functions on $\Sn$.
Obviously, 
any function
$\tilde {p}:\Sn \longrightarrow \Realsb^+$, where $\Realsb^+=[0,\infty],$
defines a 
family of confocal paraboloids $\{P(y), ~y \in \Sn\}$
with corresponding focal parameters $\tilde {p}(y)$. Even if $\tilde {p}$
is defined only on a subset of $\Sn$ we extend it to the entire $\Sn$
by setting it $=\infty$ outside of that subset. When a convex reflector is
generated by a family of confocal paraboloids described by a function on $\Sn$
we say that the reflector is generated by that function.

It is clear that any positive function $\tilde {p}(y), ~y \in \Sn$,
such that $\tilde {p}(y_1) < \infty$ and $\tilde {p}(y_2) < \infty$ 
for at least two points, generates a reflector in $\refl$. 

\begin{definition} \label{supporting}
Let $A$ be an arbitrary set in $\Rn$. A paraboloid $P(y)$
is called {\it supporting to $A$} if 
\[A\subset B(y) ~\mbox{and}~ P(y)\bigcap A \neq \emptyset.\]
\end{definition}
The following lemma is obvious.
\begin{lemma}
Let $R \in \refl$ be a reflector defined by a function $\tilde {p}$. 
Then at every point of 
$R$ there exists at least one supporting paraboloid from the family
of paraboloids defined by function $\tilde{p}$.
\end{lemma}
\begin{lemma}\label{supp-parab}\cite{refl_geom1} Let $R \in \refl$. 
For any $y \in \Sn$ there exists a unique paraboloid $P(y)$ 
(not necessarily from the family defining $R$) supporting to $R$.
\end{lemma}
For a given reflector $R$ the family of its 
supporting paraboloids $P(y),~ y \in \Sn,$ defines $R$ uniquely. However,
not every paraboloid from the family defined by a function $\tilde{p}$
is necessarily supporting to $R$. For that reason we introduce 
the following
\begin{definition}
Let $R\in \refl$. The function
$p(y), ~y \in \Sn,$ giving the value of the focal parameter for
each $P(y)$ which is supporting to $R$ is called the 
focal function of the reflector
$R$. 
\end{definition}
For $R\in \refl$ we denote by $B(R)$ 
the compact convex body bounded by $R$.
Since the origin $\cal{O}$ is an interior point of $B(R)$, 
any reflector in $\refl$ is star-shaped relative to $\cal{O}$.
Because 
$R$ is also compact, the focal function satisfies the inequalities $0 < p < \infty$ on $\Sn$.

Using the star-shapedness of $R$ 
with respect to $\cal{O}$ we can describe
$R$ as a graph over $\Sn$ by setting 
\begin{equation} \label{polarity1}
\p (x) = sup \{\lambda \geq 0~|~\lambda x \in B(R), ~x \in \Sn\}.
\end{equation}
The function $\p$ is called the {\it radial function} (cf. \cite{Schneider}).
Then any point on $R$
is given by 
\begin{equation}\label{radial}
r(x) = \p (x)x,~ x \in \Sn.
\end{equation}
Here $x$ is treated as a point on $\Sn$ and as a unit vector in $\Rn$.
It follows at once from the definition of $R$ that
\begin{equation}\label{polarity2}
\p (x) = \inf_{y \in \Sn} \frac{p(y)}{1-\xy}.
\end{equation}

It follows from (\ref{polarity2}) that
for any given $y \in \Sn$ and all $x \in \Sn$
\[p(y) \geq \p(x)(1-\xy).\]
Since for each $y\in \Sn$ there exists a paraboloid $P(y)$ supporting
to $R$, the equality is attained at some $x\in \Sn$ and therefore,
\begin{equation}\label{polarity3}
p (y) = \sup_{x \in \Sn} \p(x)(1-\xy),~y \in \Sn.
\end{equation}

The {\it reflector map} $\g$ defined by a reflector $R$ is a 
possibly multivalued map  
$\g: \Sn \rightarrow \Sn$ such that 
\begin{equation}\label{reflector_map1}
\g( x) = \bigcup_{\{P_x(y)\}} \{y\},
\end{equation}
 where $y$ is the axis of
a paraboloid $P_x(y)$ supporting to $R$ at $r(x)$ and the union is taken over
all such paraboloids. Note that by Lemma \ref{supp-parab}
the map $\g:\Sn\rightarrow \Sn$ is "onto" for any $R \in \cal{R}$. 

It follows from (\ref{polarity2}) and (\ref{polarity3}) that
the reflector map can be defined alternatively as
\begin{equation} \label{reflector_map2}
\g( x) = \{y \in \Sn~|~ p(y) = \p(x)(1-\xy)\},~x \in \Sn.
\end{equation}

If the reflector $R$ is a $C^1$
hypersurface and $u$ its outward unit normal field then
\begin{equation}\label{snell}
y = \g(x) = x - 2 \la x,u(x)\ra u(x),
\end{equation}
which is the law of reflection. In this case, the supporting
paraboloid at $r(x)=\p(x)x$ is unique and it is 
tangent to the reflector $R$ at $r(x)$.

When $R$ is not smooth, (\ref{snell}) is still
satisfied at the point  $r(x) \in R$ 
with any of the supporting paraboloids
at $r(x)$. The normal $u$ in this case is the
normal to a supporting paraboloid at $r(x)$. Here, and everywhere below,
it is assumed that the normals to any of the supporting paraboloids
are directed outward relative to the convex sets bounded by 
these paraboloids. 

%
In order to fix the terminology recall that a support function $f$ of
a nonempty closed convex body $K \subset \Rn$ is defined by
\[f(U)=\sup_{X\in K}\{\la X,U\ra\},~~U \in \Rn.\]
Since $f$ is positively homogeneous, that is,
\[f(\lambda U)= \lambda f(U)~\forall \lambda \geq 0~\mbox{and} 
~\forall U\in \Rn,\]
it is completely defined by its values on $\Sn$. In the following,
to indicate that we are considering the restriction of the support
function to $\Sn$ we denote its argument by a small letter and refer
to it as the support function of the convex hypersurface $\partial K$.
\begin{theorem}\label{smoothness of h}The support function $h$ of a reflector $R\in \refl$ is positive and
of class $C^1(\Sn)$. The reflector $R$ can be represented as
\begin{equation} \label{param of R}
X(u)= h(u)u + \nabla h(u), ~u \in \Sn,
\end{equation}
where $\nabla$ denotes the gradient in the standard metric of $\Sn$.
\end{theorem}
To prove this theorem we will need the following 
\begin{lemma} \label{support}Let $R\in \refl$, $\p$ its radial function,
  and 
$h(u),~u\in \Sn,$ its support function.
For a fixed $u\in \Sn$ consider the 
supporting hyperplane to $R$
\[\a(u) = \{Z\in \Rn ~|~\la Z,u\ra = h(u)\}\]
and let  $C_u=R\bigcap \a(u)$. 
Then there exists
only one $x \in \Sn$ such that
\begin{equation} \label{!linseg}
C_u = \{\p(x)x\}.
\end{equation}
Furthermore, there exists no more than one paraboloid $P$
supporting to $R$ for which $\a(u)$ is the tangent hyperplane at
$\{\p(x)x\}$.
\end{lemma}
\begin{proof} Fix some $u \in \Sn$ and let $X \in C_u$. 
Let $P_X$ be a paraboloid supporting to $R$ at $X$ and $TP_X$ the
tangent hyperplane to $P_X$ at $X$. Since $TP_X\bigcap P_X = \{X\}$,
we have $TP_X \bigcap R = \{X\}$, that is, $TP_X$ is also a supporting
hyperplane to $R$ at $X$. Thus, any hyperplane tangent to a
supporting paraboloid may contain at most one point of
$R$. If 
$TP_X=\a(u)$, then (\ref{!linseg}) is clearly true. 

Suppose, $TP_X\neq\a(u)$ and let $X_1 \in C_u,~X_1\neq X$. Since $R$ is
a boundary of a compact convex body, the segment
$XX_1: X_t = (1-t)X+tX_1,~t \in [0,1]$, is contained in $C_u$. Clearly,
for any
$t \in (0,1)$ there exists no paraboloid supporting to $R$ at  $X_t$,
because for any such paraboloid its tangent hyperplane at $X_t$ will
contain the segment $XX_1$, which is impossible. Thus, 
$R\bigcap \a(u)=\{X\}$ and this implies (\ref{!linseg}). 

The last statement of the lemma follows essentially from the preceding part. 
Indeed, for any paraboloid $P$ supporting to $R$ at $X$ and  
such that the tangent hyperplane at
the point of support is $\a(u)$, the contact set 
$R\bigcap \a(u)=\{X\}$. Let $x =X/|X|$.
By reflection law (\ref{snell}) the axis of $P$ is
defined uniquely as
\[y = x -2\la x,u\ra u.\]
Since the focus of $P$ is fixed at $\cal O$, it follows from
(\ref{polrad}) that the vectors $x,y$ define $P$ up to a homothety
with respect to $O$. But $P$ also contains $X$, and therefore it
is defined uniquely. 
\end{proof}

\begin{proofof}$~${\bf of Theorem \ref{smoothness of h}.} 
By Lemma \ref{support},
 for any $u \in \Sn$ the contact set $C_u=R\bigcap \a(u)$
consists of only one point. Now, except for positivity of $h$, all statements
 of this theorem follow
from Corollary 1.7.3 in \cite{Schneider}, p. 40.

The positivity of $h$
follows from the fact that the origin $\calO$ is strictly inside the
convex body bounded by $R$. 
\end{proofof}


\section{\protect\bf A Minkowski-type theorem for focal functions}
\label{mink}
The main goal of this section is to determine conditions under which a given function
on $\Sn$ is the focal function of a convex reflector in $\refl$.
The following theorem answers this question.
\begin{theorem} \label{minkg} Let $R\in \refl$ and let 
$p$ be the focal function of $R$. Then $p$ is a positive bounded function
on $\Sn$ and if for  $x, y\in \Sn$ 
\begin{equation}\label{minkg1}
x-y = \sum_{i=1}^k\a_i(x-y_i),~~\a_i \geq 0,~~\sum_{i=1}^k\a_i > 0,
\end{equation}
where the vectors $\{x-y_1,...,x-y_k\}$ are linearly independent,
then the inequality
\begin{equation} \label{minkg2}
p(y) \leq \sum_{i=1}^k\a_ip(y_i)
\end{equation}
is satisfied. 

Conversely, let $p$ be a positive bounded function
on $\Sn$ such that for  $x, y\in \Sn$
satisfying (\ref{minkg1}) the inequality (\ref{minkg2}) holds.
Then there exists a unique reflector $R\in \refl$
with the focal function $p$.
\end{theorem} 
We shall need the following 
\begin{lemma}\label{decomp}
Let $R \in \refl$, $y_0 \in \Sn,$ and $P(y_0)$ a paraboloid supporting
to $R$ at some point $X$. Let $x= X/|X|$. Then the following
decomposition holds:
\begin{equation}\label{decomp0}
x - y_0 = \sum_{i=1}^k \a_i(x-y_i), ~~\a_i \geq 0,~~\sum_{i=1}^k\a_i > 0,,~~\mbox{and}~ k \leq n+1,
\end{equation}
where $y_i$ are the axes of paraboloids
supporting to $R$ at $X$,
vectors $x-y_1, ...,x-y_k$ are linearly independent and $k \leq n+1$.
\end{lemma}
\begin{remark} The case when $y_0$ coincides with one of the
vectors $y_1,...,y_k$ is not excluded. In such case there is
only one term in the sum and the corresponding coefficient is
equal to one.
\end{remark}
\begin{proof} Consider the set 
of rays originating at $X$ and
intersecting $R$ at points different from $X$. The closure of this
set (as a set of directions on a unit sphere centered at $X$) 
we denote by $\Psi_X$ and its boundary by $\psi_X$.
 Because $\calO$ is an interior point of $B(R)$, 
the cone $\Psi_X$
has a nonempty interior. It is also convex. 
The rays in $\psi_X$ form the tangent cone to $R$ at $X$.


Let $\Phi_X$
be the cone with vertex $X$ dual to $\Psi_X$. We denote
its boundary by $\phi_X$. It is a cone dual to $\psi_X$. Because 
${\Psi}_X$ has a nonempty interior and $\Phi_X$ is convex, it is clear that
there exists a hyperplane $Q$ passing through $X$ for which
$Q\bigcap \Phi_X = \{X\}$. Let $Q^+$ be the halfspace
determined by $Q$ which contains $\Phi_X$ and
$Q_1$ the hyperplane
parallel to $Q$, at a distance $1$ from $Q$, in the halfspace $Q^+$.
The intersection $\Pi = Q_1\bigcap\Phi_X$ is
a compact convex set on $Q_1$. Its boundary $\pi=Q_1\bigcap\phi_X$.

Observe also, that because $\cal O$ is in the
interior of $B(R)$, Theorem \ref{smoothness of h} implies that for
all hyperplanes supporting to $R$ at $X$
we have
\begin{equation}\label{minimum h}
\la x,u\ra = \frac{h(u)}{|X|} \geq c_0 =const > 0,
\end{equation}
where $u$ is the outward normal to such
a hyperplane and $c_0$ is independent
on a particular supporting hyperplane at $X$.

Let now $P(y_0)$ be a supporting paraboloid at $X$. It follows
from (\ref{snell}) applied to $P(y_0)$ that $x-y_0$ is orthogonal to
the hyperplane tangent to $P(y_0)$ at $X$ and therefore the ray originating
at $X$ in direction $x-y_0$ is in the set $\Phi_X$. Let
$v_0$ denote the point of intersection of this ray with $Q_1$.
Since $\Pi$ is the  convex hull of $\pi$ and $v_0 \in \Pi$,
by Carath\'{e}odory's theorem we have,
\begin{equation}\label{carat}
v_0 = \sum_{i=1}^k\beta_iv_i, ~\beta_i \geq 0, \sum_{i=1}^k\beta_i =1,
~v_i \in \pi,
\end{equation}
where $v_1,...,v_k$ are affinely independent as points of $Q_1$ 
(considered as a Euclidean space)  and $k \leq n+1$.

On the other hand, each $v_i \in \pi$ is a point on the ray
orthogonal to some hyperplane containing a ray tangent to $R$ at
$X$. Such hyperplane is also tangent to some paraboloid $P(y_i)$
supporting to $R$ at $X$. Therefore,
for some $\gamma_i > 0$
\[v_i = X + \gamma_i (x - y_i)=X + 2\gamma_i\la x,u_i\ra u_i,\]
where $u_i$ is the outward unit normal to the
hyperplane tangent to $P(y_i)$ at $X$. In addition, 
by (\ref{minimum h}), $\gamma_i < \infty$. Similarly,
for some $\gamma_0, ~0 < \gamma_0 < \infty$, $v_0 = X + \gamma_0(x-y_0)$.
This together with (\ref{carat}) implies (\ref{decomp0}) with
$\a_i = \beta_i\gamma_i/\gamma_0 \geq 0$. It also follows from (\ref{carat}) that
\[\sum_{i=1}^k \a_i > 0.\]

Finally, the linear independence of vectors 
$\gamma_1 (x - y_1),..., \gamma_k (x - y_k)$ in $\Rn$ follows from affine independence of points $v_1,...,v_k$ in $Q_1$.
\end{proof}

\begin{proofof}$~${\bf of Theorem \ref{minkg}. Necessity.} 
Let $R\in \refl$ and let $p(y), ~y \in \Sn,$ be the focal function of
$R$. The boundedness of $p$ follows from compactness of $R$
and $p >0 $ on $\Sn$ because $\calO$ is an interior point of $B(R)$
and by (\ref{polarity3}).

We show now that the inequalities (\ref{minkg2}) are satisfied.
Indeed, let $X$ be an arbitrary point on $R$,
$x= X/|X|$ and $P(y)$ a paraboloid supporting to $R$ at $X$.
Then, taking into account (\ref{minkg1}), we have
\[p(y) = |X|\la x-y,x\ra=\sum_i\alpha_i (|X|\la x-y_i,x\ra).\]
By Lemma \ref{supp-parab} for each $y_i$ there exists a
supporting paraboloid $P(y_i)$. By definition
of the focal function the  focal parameters of $P(y_i)$ are
 $p(y_i)$ and, since $P(y)$ is supporting at $X$,
we have by (\ref{polarity2})
\[|X| \leq \frac{p(y_i)}{\la x-y_i,x\ra}.\]
Then
\[p(y)\leq\sum_i\alpha_i p(y_i).\]
This completes the proof of necessity.

\noindent {\bf Sufficiency.} Since $p$ is positive and bounded,
we can define the reflector $R$ as
\[R = \partial (\bigcap_{y\in \Sn}B(y)),\]
where $B(y)$ is the closed convex body bounded by paraboloid
$P(y)$ with focal parameter $p(y)$. Obviously, $R\in \refl$.
We need to show that for each $y\in \Sn$ $p(y)$ is a focal
parameter of a paraboloid supporting to $R$. 

Fix any such $y$ and suppose the corresponding paraboloid is
not supporting. By construction, $R \subset B(y)$ and, 
since $P(y)$ is not supporting to $R$, $R \bigcap P(y)=\emptyset$.
By Lemma \ref{supp-parab} there exists a unique paraboloid ${\bar P}(y)$
supporting to $R$ at some $X\in R$. Let ${\bar p}(y)$ be the focal
parameter of ${\bar P}(y)$. Obviously, ${\bar p}(y) < p(y)$.

By Lemma \ref{decomp} we have the decomposition
\[
x-y = \sum_{i=1}^k \a_i(x-y_i), ~~\a_i \geq 0,~~\sum_{i=1}^k\a_i > 0,~\mbox{and}~ k \leq n+1,
\]
where $x=X/|X|$ and $y_i$ are the axes of paraboloids
supporting to $R$ at $X$ and not all of $\a_i$ are zeroes. Then
\[
{\bar p}(y)= |X|\la x-y,x \ra = \sum_{i=1}^k \a_i|X|\la x-y_i,x\ra
= \sum_{i=1}^k \a_ip(y_i) < p(y),
\]
which contradicts (\ref{minkg2}).
\end{proofof}

It is of interest to compare Lemma \ref{decomp} 
with the well known result for convex hypersurfaces which states that
that if $M$ is a convex hypersurface, 
$Q$ is a supporting hyperplane to $M$ at $X$ and 
$u$ is the  outward unit normal to $Q$
then there exist $\a_i, i=1,...,k \leq n+1, ~\a_i \geq 0,$ such that the decomposition
\[ u = \sum_{i=1}^{k}\a_iu_i\]
holds with $u_1,...,u_k$ being linearly independent outward normals to 
hyperplanes supporting to $M$ at $X$ \cite{Bak}, p. 40. The seemingly natural
extension of this result to convex reflectors with normals replaced
by the axes of supporting paraboloids is not possible. This can be
seen from the following example.

Let $P_1=P(y_1), ~P_2=P(y_2)$ be two confocal paraboloids
 in $\Reals^3$, with axes $y_1$ and $y_2$, respectively, common focus
$\calO$ and positive focal parameters. Assume that $y_1$ and $y_2$ are perpendicular
to each other. Let $R_2$ be the convex reflector formed by
$P_1$ and $P_2$. Consider a point $X \in P_1\bigcap P_2$ not
lying in the plane spanned by $y_1$ and $y_2$. The axes of
paraboloids which are supporting to $R_2$ at $X$
lie on a circular cone with axis parallel to the
line tangent to $P_1\bigcap P_2$ at $X$. Any such axis
different from $y_1$ and $y_2$ can not be represented
as a linear combination of $y_1$ and $y_2$.

In the rest of this section we consider an extension of a focal function
to the entire $\Rn$ and give a quantitave expression for its
``deviation'' from a subadditive function.

Recall that a function $f: \Rn \longrightarrow \Reals$ is 
$sublinear$ if it is positively homogeneous
and $subadditive$. The latter means that 
\[f(X + X^{\prime}) \leq f(X) + f(X^{\prime})~\forall X, X^{\prime} \in \Rn.\]
A well known theorem of Minkowski asserts that sublinearity
is a necessary and sufficient condition for a function in
$\Rn$ to be the
support function of a unique compact convex body \cite{Schneider}, p. 38.

Let  $R \in \refl$ and let $p$ be its focal function.
It follows from the definition (\ref{polarity1}) of the radial function
$\p$ of $R$ and from
the expression (\ref{polarity3}) for its focal function $p$,
that $p$ can be naturally extended to the
entire $\Rn$ by setting
\begin{equation}\label{ext2Rn}
p(Y) = \sup_{X \in B(R)} (|X||Y| - \XY ), ~Y \in \Rn.
\end{equation}
The extended $p$ is obviously positively homogeneous and
it is shown in \cite{refl_geom1} that it
 is also subadditive. Then, by Minkowski's theorem, it is the support function
of a compact convex body in $\Rn$. It will be shown 
in the next section (Proposition
\ref{param for D(R)}) that the boundary of this body is
the directrix hypersurface (defined in the next section) of the reflector $R$.

The next proposition shows that focal functions of convex reflectors
satisfy an inequality more restrictive than the one defining subadditivity. 
\begin{proposition} \label{subadd0}Let $R\in \refl$ and let 
$p(y),~y\in \Sn$,
be its focal function. Extend $p$ to the entire $\Rn$ as
in (\ref{ext2Rn}). Then for any $Y_1,...,Y_N \in \Rn$ and $\a_i \geq 0,~i=1,...,N,$
the extended focal function of $R$ satisfies
the inequality
\begin{equation}\label{jenseng}
p(\sum_{i=1}^N\a_iY_i) \leq \sum_{i=1}^N\a_ip(Y_i) +
|X|(|\sum_{i=1}^N\a_iY_i|-\sum_{i=1}^N\a_i|Y_i|),
\end{equation}
where $X$ is a point on $R$ 
at which the paraboloid with axis $\sum_{i=1}^N\a_iY_i$ is
supporting to $R$. It is assumed here that if 
$\sum_{i=1}^N\a_iY_i=\cal O$
then $X=\cal O$.
\end{proposition}
\begin{proof} 
First, note that the right hand side in (\ref{jenseng}) is always
nonnegative. Indeed, using (\ref{ext2Rn}), we obtain
\begin{eqnarray}
\sum_{i=1}^N\a_ip(Y_i) +
|X|(|\sum_{i=1}^N\a_iY_i|-\sum_{i=1}^N\a_i|Y_i|) \nonumber \\
=\sum_{i=1}^N\a_i\sup_{Z\in B(R)}(|Z||Y_i|-\la Z,Y_i\ra) +
|X|(|\sum_{i=1}^N\a_iY_i|-\sum_{i=1}^N\a_i|Y_i|) \nonumber \\
\geq  
\sum_{i=1}^N\a_i(|X||Y_i|-\la X,Y_i\ra) +
|X|(|\sum_{i=1}^N\a_iY_i|-\sum_{i=1}^N\a_i|Y_i|) \nonumber \\
=|X||\sum_{i=1}^N\a_iY_i|-\la X,\sum_{i=1}^N\a_iY_i \ra \geq 0. \nonumber
\end{eqnarray}

Put
$Y = \sum_{i=1}^N\a_iY_i.$
Then for the paraboloid with axis $Y$ supporting to $R$ at $X$, we have
\begin{eqnarray}
p(Y) = |X||Y|-\la X, Y\ra = \sum_{i=1}^N\a_i(|X||Y_i|-\la X, Y_i\ra) +
|X|(|Y|-\sum_{i=1}^N\a_i|Y_i|)\nonumber \\
\leq \sum_{i=1}^N\a_ip(Y_i) + |X|(|Y|-\sum_{i=1}^N\a_i|Y_i|).\nonumber
\end{eqnarray}
\end{proof}
\begin{remark} \label{subadd1}Taking $Y= Y_1 +Y_2$ in the above proposition one sees immediately
that the extended focal function is subadditive. The proposition also shows that
subadditivity by itself is not sufficient for a positive and bounded function to
be a focal function of a convex reflector.
\end{remark}
\begin{remark} If $Y= \sum_{i=1}^N\a_iY_i,~\a_i \leq 0,~i=1,2,...,N,$
then
\begin{equation}\label{jenseng1}
p(\sum_{i=1}^N\a_iY_i) \geq \sum_{i=1}^N\a_ip(Y_i) +
|X|(|\sum_{i=1}^N\a_iY_i|-\sum_{i=1}^N\a_i|Y_i|).
\end{equation}
The proof is similar to the proof of Proposition \ref{subadd0}.
\end{remark}
\section{\protect\bf Directrix of a reflector}
\label{dir}
Associated with a reflector $R\in \cal{R}$ is its {\it directrix} hypersurface
$D(R)$
defined as follows. 
Let $r(x)=\p(x)x, ~x \in \Sn,$ be the position vector
of the reflector $R$, where $\p$ is the radial function. For 
each $x \in \Sn$ put
\[D_{x}(R) = r(x) - \bigcup_{y\in \g(x)} \{\p(x)y\},\]
where $\g$ is the reflector map. The directrix
of $R$ is defined as
\[D(R)=\bigcup_{x \in \Sn}D_{x}(R).\]
The directrix of a reflector was introduced and partially studied in
\cite{refl_geom1}. The Theorem \ref{param for D(R)} below
provides additional and more detailed information about $D(R)$.
\begin{theorem}\label{param for D(R)} Let $R \in \refl$. Then:
\begin{enumerate}
\item $D(R)$ 
is a closed convex hypersurface with the origin $\cal O$  in
the interior of the compact convex body bounded by $D(R)$;
\item for any paraboloid
$P(y)$ supporting to $R$
the vector $-y$ is the outward unit normal
to a hyperplane supporting to $D(R)$.
The hyperplanes supporting to $D(R)$
are the directrix hyperplanes of paraboloids supporting to $R$;
\item the support function $H$ of $D(R)$
satisfies the equality
\begin{equation}
H(-y) = p(y),~y \in S^n,
\end{equation} where $p$ is the focal function of $R$;
\item $D(R)$ can be parametrized as
\begin{equation}\label{param for D(R)_1}
r^D(u)=2h(u)u,~u \in S^n,
\end{equation}
where $h$ is the support function of $R$.
In this parametrization $D(R)$ is of class $C^1$.
\end{enumerate}
\end{theorem}
\begin{proof} Remark \ref{subadd1}  together  with
the positive homogeneity of the focal function extended
to $\Rn$ imply that $p$ is sublinear and, therefore, it is
the support function of a unique compact convex body in $\Rn$. 
Denote this body by $\cal D$. Since $p > 0$,
the origin $\cal O$ is an interior point of $\cal D$ and, thus,
$\cal D$ is nondegenerate. 

We show now that the boundary of $\cal D$ is the directrix $D(R)$. 
Fix some $\bar{y} \in \Sn$ and consider the hyperplane 
\[Q(\bar{y}) = \{Z \in \Rn | \la Z, -\bar{y}\ra = p(\bar{y})\}\]
and the halfspace
\[Q^{-}(\bar{y})= \{Z \in \Rn | \la Z, -\bar{y}\ra 
 \leq
 p(\bar{y})\}.\]
We want to show that $Q(\bar{y})$ is supporting
to $D(R)$.

Let $P(\bar{y})$
be  the paraboloid with axis $\bar{y}$ 
supporting to $R$.  By Lemma
\ref{supp-parab} such paraboloid exists and unique. 
Then
\[|X|\la x - \bar{y}, x\ra) \leq p(\bar{y})~~\forall X \in B(R)\]
and the equality is attained only when $X  \in R \bigcap P(\bar{y})$; 
here, as usual, $x=X/|X|$ for $X \neq \calO$. When $X = \calO$
this inequality, with any $x \in \Sn$, is obvious as $p$ is positive.
Since any $Z \in D(R)$ is given by
\[Z = |X|(x - y)\]
for some $X\in R$ with $y$ being the axis of a supporting paraboloid at
$X$, we have
\[\la Z, -\bar{y}\ra = |X|\la x - y,-\bar{y}\ra \leq\la X,-\bar{y}\ra+|X|
\leq p(\bar{y}),\]
that is, $D(R) \subset Q^{-}(\bar{y})$.

Let $\bar{X} \in R \bigcap P(\bar{y})$. 
Then $|\bar{X}|(\bar{x} - \bar{y}) \in D(R)\bigcap Q(\bar{y})$,
that is, the hyperplane $Q(\bar{y})$ is supporting to $D(R)$.
Thus, for each $\bar{y} \in \Sn$ 
$D(R)$ and $\cal D$ have the same supporting hyperplane 
$Q(\bar{y})$.
Hence, $\partial {\cal D} = D(R)$ and
the support function $H(-y)$ of $D(R)$ is $p(y)$. It is also
clear that $Q(\bar{y})$ is the directrix hyperplane of $P(y)$. This
proves 1-3.

Let us show now that $D(R)$ can be parametrized as in 
(\ref{param for D(R)_1}). Let $u \in \Sn$ and let
$\a(u)$ be the hyperplane
with outward unit normal $u$ supporting to $R$. By Lemma \ref{support}
there exists a unique $X(u) \in R$ 
such that $\a(u)\bigcap R = \{X(u)\}$. Consider a ray of direction
$u$ originating at $\cal O$. Since $D(R)$ is convex and 
$\cal O$ is in the interior of the compact convex body bounded by
$D(R)$, this ray intersects $D(R)$ in a unique point $Z(u)$ and
by definition of $D(R)$ 
\[Z(u) = X(u) - |X(u)|y\]
for some $y \in \g(x(u))$. But
\[u = \frac{Z(u)}{|Z(u)|} = \frac{x(u)-y}{|x(u)-y|}.\]
By the reflection law (\ref{snell}), applied to paraboloid
$P(y)$ supporting to $R$ at $X(u)$,  $u$ is the normal
vector to $P(y)$ at $X(u)$, that is, the hyperplane $\a(u)$
is also the tangent hyperplane to $P(y)$ at $X(u)$.
Then, 
\[\la Z(u),u\ra = |X(u)|\la x(u)-y,u\ra = 2|X(u)|\la x(u),u\ra = 2h(u)u.\]
This proves (\ref{param for D(R)_1}).

The $C^1-$smoothness of $D(R)$ follows now from Theorem 
\ref{smoothness of h}.
\end{proof}
\newpage
\begin{remark}
The representation (\ref{param for D(R)_1}) shows that
the directrix $D(R)$ is a rescaled (with coefficient 2) 
pedal hypersurface of $R$; cf. \cite{Hasanis_Koutrou:85} and 
other references there.
\end{remark}
\bibliographystyle{plain} 
\bibliography{oliker.bbl}

\begin{thebibliography}{10}

\bibitem{Bak}
I.~J. Bakelman.
\newblock {\em Convex Analysis and Nonlinear Geometric Elliptic Equations}.
\newblock Springer-Verlag, Berlin, 1994.

\bibitem{CGH:04}
L.~A. Caffarelli, C.~Gutierrez, and Qingbo Huang.
\newblock On the regularity of reflector antennas.
\newblock {\em Preprint}, 2004.

\bibitem{Caf_alloc:96}
L.A. Caffarelli.
\newblock Allocation maps with general cost functions.
\newblock {\em Partial Differential Equations and Applications}, 177:29--35,
  1996.

\bibitem{CKO}
L.A. Caffarelli, S.~Kochengin, and V.I. Oliker.
\newblock On the numerical solution of the problem of reflector design with
  given far-field scattering data.
\newblock {\em Contemporary Mathematics}, 226:13--32, 1999.

\bibitem{CO}
L.A. Caffarelli and V.I. Oliker.
\newblock Weak solutions of one inverse problem in geometric optics.
\newblock {\em Unpublished manuscript}, 1994.

\bibitem{galo:04}
W. Gangbo and V. I. Oliker.
\newblock Existence of optimal maps in the reflector-type problems.
\newblock {\em ESAIM: Control, Optimization and Calculus of Variations}, 
to appear.

\bibitem{glimm-oliker-refl2:03}
T.~Glimm and V.I.~Oliker.
\newblock Optical design of single reflector systems and the
  {M}onge-{K}antorovich mass transfer problem.
\newblock {\em J. of Math. Sciences}, 117(3):4096--4108, 2003.

\bibitem{glimm-oliker-refl1:03}
T.~Glimm and V.I.~Oliker.
\newblock Optical design of two-reflector systems, the {M}onge-{K}antorovich
  mass transfer problem and {F}ermat's principle.
\newblock {\em Indiana Univ. Math. J.}, 53:1255--1278, 2004.

\bibitem{Hasanis_Koutrou:85}
T.~Hasanis and D.~Koutroufiotis.
\newblock The characteristic mapping of a reflector.
\newblock {\em J. of Geometry}, 24:131--167, 1985.

\bibitem{KA}
L.V. Kantorovich and G.P. Akilov.
\newblock {\em Functional Analysis, ch. VIII, \S 4, In Russian}.
\newblock Nauka, Moscow, 1977, 2-nd revised edition.

\bibitem{refl_geom1}
V.I. Oliker.
\newblock On the geometry of convex reflectors.
\newblock {\em PDE's, Submanifolds and Affine Differential Geometry, ed. by B. Opozda, U. Simon and M. Wiehe, Banach
  Center Publications}, 57:155--169, 2002.
\newblock {Errata: Banach Center Publications, v. 69(2005), 269-270}.

\bibitem{patow_pueyo}
G.~Patow and X.~Pueyo.
\newblock A survey of inverse surface design from light transport 
behavior specification.
\newblock {\em Computer Graphics Forum}, 24:773--789, 2005.

\bibitem{Schneider}
R.~Schneider.
\newblock {\em {Convex Bodies. The Brunn-Minkowski Theory}}.
\newblock Cambridge Univ. Press, Cambridge, 1993.

\bibitem{Wang:inv03}
Xu-Jia Wang.
\newblock {On design of a reflector antenna II}.
\newblock {\em Calculus of Variations and PDE's}, 20:329--341, 2004.

\bibitem{W}
B.~S. Westcott.
\newblock {\em Shaped Reflector Antenna Design}.
\newblock Research Studies Press, Letchworth, UK, 1983.

\end{thebibliography}
\end{document}